\documentclass[12pt,reqno]{amsart}

\usepackage{amssymb,latexsym}
\usepackage[matrix,arrow]{xy}

\theoremstyle{plain} 
\newtheorem{thm}{Theorem}[section]
\newtheorem{cor}[thm]{Corollary}

\newtheorem{prop}[thm]{Proposition}
\theoremstyle{remark} 
\newtheorem{rem}[thm]{Remark}
\theoremstyle{definition}

\newcommand{\rl}{\mathbb{R}}
\newcommand{\cx}{\mathbb{C}}
\newcommand{\cl}{\mathbb{C} l}
\newcommand{\set}[1]{\left\{#1\right\}}
\newcommand{\so}{\mathrm{SO}}
\newcommand{\spin}{\mathrm{Spin}}

\newcommand{\wit}[1]{\widetilde{#1}}
\newcommand{\ga}{\alpha}

\newcommand{\gs}{\sigma}
\newcommand{\go}{\omega}
\newcommand{\gl}{\lambda}
\newcommand{\gS}{\Sigma}
\newcommand{\gp}{\varphi}
\newcommand{\gn}{\nabla}

\newcommand{\ovl}[1]{\overline{#1}}

\newcommand{\demi}{\frac{1}{2}}
\newcommand{\gG}{\Gamma}

\newcommand{\quart}{\frac{1}{4}}
\newcommand{\ls}{\setlength{\baselineskip}{16pt}
                 \setlength{\parskip}{3mm}}

\title[]{\bf Eigenvalue Estimates for the Dirac-Schr\"{o}dinger Operators}

\author[]{Bertrand Morel} 
\address{Institut {\'E}lie Cartan\\
Universit{\'e} Henri Poincar{\'e}, Nancy I\\
B.P. 239\\
 54506 Vand\oe uvre-L{\`e}s-Nancy Cedex, France}
\email{morel@iecn.u-nancy.fr}



\begin{document}
\ls 
\begin{abstract}
We give new estimates for the eigenvalues of the hypersurface Dirac operator in 
terms of the intrinsic energy-momentum tensor, the mean curvature and the scalar 
curvature. We also discuss their limiting cases as well as the limiting cases of 
the estimates obtained by X.~Zhang and O.~Hijazi in \cite{Z1} and \cite{HZ1}. We 
compare these limiting cases with those corresponding to the Friedrich and 
Hijazi inequalities. We conclude by comparing these results to intrinsic 
estimates for the Dirac-Schr\"{o}dinger operator $D_f=D-\frac{f}{2}$. \vspace{0.5cm} 

\end{abstract}
\maketitle
\section{Introduction}
In this paper, we start by comparing the hypersurface spinor bundle $S$ of a 
hypersurface $M$ to the fundamental spinor bundle $\gS M$ of $M$. The 
hypersurface spinor bundle $S$ is obtained by restricting the spinor bundle of 
the ambient space $N$ to $M$. If $\gp \in \gG(S)$ is a section of this bundle, 
the energy-momentum tensor $Q^{\,\gp}$ associated with $\gp$ is defined on the 
complement of its zero set, by $$Q^{\,\gp}_{ij}=\frac{1}{2} (e_i \cdot \nu 
\cdot\gn_j \gp + e_j \cdot\nu \cdot\gn_i \gp, \gp/|\gp |^2),$$ where $\nu$ is a 
unit normal vector field globally defined along $M$, $e_i,e_j$ are vectors of a local 
orthonormal frame of $M$, and where $\gn_i \gp$ stands 
for the covariant derivative of the spinor field $\gp$ in the direction of $e_i$ . Then the 
Schr\"odinger-Lichnerowicz formula for the classical Dirac operator $D$ on $M$ 
leads to the following result (compare with \cite{Z1}): 

\begin{thm}\label{mthm}
Let $M^n \subset (N^{n+1},\wit{g}\,)$ be a compact hypersurface of a Riemannian 
spin manifold $N$.  Let $\gl$ be any eigenvalue of the hypersurface Dirac 
operator $D _H=D-\frac{H}{2}$ , associated with an eigenspinor $\gp$. Assume 
that $R+4 |Q^{\,\gp}| ^2 
> H^2 > 0$, then one has 
\begin{equation}\label{bin}
\gl^2 \geq \frac{1}{4}\; \inf _M \left( \sqrt {R+4 |Q^{\,\gp}| ^2} -|H| \right) ^2. 
\end{equation}
where $R$ and $H$ are respectively the scalar curvature and the mean curvature 
of $M$, and $Q^{\,\gp}$ is the energy-momentum tensor associated with $\gp$. 
\end{thm}

In fact, we see that if $M$ is a minimal hypersurface, the hypersurface Dirac 
operator corresponds to the classical Dirac operator. Therefore, in this case, 
this estimate is exactly the one given by O. Hijazi in \cite{Hi3}. 

We then discuss the limiting case of equation (\ref{bin}) and that given by X. 
Zhang in \cite{Z1}.  

As in  \cite{Hi1} and \cite{HZ1}, we then prove 

\begin{thm}\label{sthm}
Under the same conditions as in Theorem \ref{mthm}, suppose that $ \overline R\; 
e ^{2u} +4|Q^{\,\gp}| ^2>H ^2>0$, where $\overline R$ is the scalar curvature of 
$M$ for some conformal metric $\overline g = e ^{2u}\,\wit{g}$, with 
$du(\nu)_{\mid M}\equiv 0$, then 
\begin{eqnarray}\label{ssbin}
\gl^2 \geq \frac{1}{4}\;  \inf _M \left( \sqrt {\overline R\; e ^{2u}+4 
|Q^{\,\gp}| ^2} -|H| \right) ^2. 
\end{eqnarray}
\end{thm}

The discussion of the limiting case in this inequality and that proved in 
\cite{HZ1} is similar to that of (\ref{bin}). As a conclusion, we observe that 
these inequalities correspond to a generalization of the classical estimates in 
terms of the Dirac-Schr\"{o}dinger operator $D_f=D-\frac{f}{2}$, for a real function 
$f$ on $M$. 

We would like to thank Oussama Hijazi for pointing out this problem, as well as 
Nicolas Ginoux and Xiao Zhang for helpful discussions. 

\section{Preliminaries}

\subsection{Restriction of Spinors to the Hypersurface}


In this paper we will consider an oriented compact hypersurface $(M^n,g)$ of a 
Riemannian spin manifold $(N^{n+1},\wit{g})$, with a spin structure $\spin N$. 
The metric $g$ is the induced metric on $M$ by $\wit{g}$. The possibility to 
define globally a unit normal vector field $\nu$ on $M$ allows to induce from 
$\spin N$ a spin structure on $M$, denoted by $\spin M$. For this, we can 
associate to every oriented orthonormal frame $(e_1,\ldots,e_n)$ on $M$ an 
oriented orthonormal frame $(e_1,\ldots,e_n,\nu)$ of $N$ such that the principal 
$\so(n)$-bundle $\so_nM$ of oriented orthonormal frames on $M$ is identified 
with a sub-bundle of $\so_{n+1} N_{\mid M}$. Such a map is denoted by $\Phi$.  

Let $\cl_n$ be the $n$-dimensional complex Clifford algebra and $\cl_n^0$ its 
even part. Recall that there exists an isomorphism 
\begin{eqnarray}\label{alpha1}
 \ga : & \cl_n &\longrightarrow  \cl_{n+1}^0 \\  
  & e_i &\longmapsto    e_i\cdot \nu.\nonumber
\end{eqnarray}
Here, $\nu$ stands for the last vector of the 
canonical basis of  $\rl^{n+1}$. 

In particular, $\ga$ yields the following commutative diagram: $$\xymatrix{\spin 
(n) \ar[d]^{\mathrm{Ad}}\ar[r]^\ga&\spin (n+1)\ar[d]^{\mathrm{Ad}}\\ 
\so(n)\ar@{^{(}->}[r]&\so (n+1)},$$ where the inclusion of $\so(n)$ in 
$\so(n+1)$ is that which fixes the last basis vector under the action of 
$\so(n+1)$ on $\rl^{n+1}$, and $\mathrm{Ad}$ the adjoint representation of 
$\spin (n)$ on $\so(n)$, which is given by $$\mathrm{Ad}_\eta(x)=\eta\cdot 
x\cdot\eta^{-1}$$ for all $\eta \in \spin (n)$ and $x \in \rl^n$. 

This allows to pull back via $\Phi$ the fiber bundle $\spin N_{\vert M}$ on 
$\so M$ as a spin structure for $M$, denoted by $\spin M$. The projection of 
$\spin M$ on $\so M$, as well as the projection of $\spin N$ 
on $\so N$, is denoted is $\pi$. Thus, we have he following commutative diagram:$$\xymatrix{\spin M 
\ar[d]^\pi\ar[r]^{\Phi^*}&\spin N_{\mid M}\ar[d]^\pi\\ \so 
M\ar@{^{(}->}[r]^\Phi&\so N_{\mid M}},$$ 

Let $\gS N$ be the spinor bundle on $N$, i.e., $$\gS N=\spin N 
\times_{\rho_{n+1}}\gS_{n+1},$$ where $\rho_{n+1}$ is the restriction to 
$\spin(n+1)$ of an irreducible representation of $\cl_{n+1}$ on the space of 
spinors $\gS_{n+1}$, of dimension $2^{[\frac{n+1}{2}]}$ ($[\,.\,]$ denotes the 
integer part). Recall that if $n+1$ is odd, this representation is chosen so 
that the complex volume form acts as the identity on $\gS_{n+1}$. 

Locally, by definition of $\gS N$, if $U$ is an open subset of $N$ and $\psi \in 
\gG_U(\gS N)$ a local section of the spinor bundle, we can write 
$$\psi=[\wit{s},\gs]$$ where $\gs : U \rightarrow \gS_{n+1}$ and $\wit{s} : U 
\rightarrow \spin N$ are smooth maps, and $[\wit{s},\gs]$ is the equivalence 
class with respect to the relation 
$$[\wit{s},\gs]\sim[\wit{s}g,\rho_{n+1}(g^{-1})\gs], \qquad\forall g\in\spin 
(n+1).$$ Moreover, we can always assume that $\pi(\wit{s})$ is a local section 
of $\so N$ with $\nu$ for last basis vector. Then we have $$\psi_{\vert 
M}=[(\wit{s}_{\vert U\cap M},\gs_{\vert U\cap M})]$$ where the equivalence class 
is reduced to elements of $\spin (n)$. 

It follows that one can realize the restriction to $M$ of the spinor bundle $\gS 
N$ as $$S:=\gS N_{\vert M}=\spin M \times_{\rho_{n+1}\circ\ga}\gS_{n+1}.$$ 

{\bf Remark.} The inclusion of $\spin(n)$ in $\spin(n+1)$ given by $\ga$ is the 
trivial one. But, this notation emphasizes that Clifford multiplication of a 
spinor field $\phi\in \gG (S)$ by a vector $X$ tangent to $M$ is given by 
\begin{equation}\label{mulcli}
(X,\phi)\mapsto X\cdot \nu\cdot\phi. 
\end{equation}
This fact is crucial for the following identification (see \cite{Bar1}, 
\cite{BHMM}). 

\subsection{Identification of $S$ with $\gS M$}

We now compare $S$ with the intrinsic spinor bundle of $M$, $$\gS M= \spin 
M\times_{\rho_n}\gS_n.$$ For this, we have to examine the cases where $n$ is 
even or odd. First assume that $n=2m$ is even. From (\ref{alpha1}) and
\begin{equation}\label{cl2m}
  \cl_{2m}\cong\cx(2^m),
\end{equation}
it follows that the representation of $\cl_{2m}$ given by $\rho_{2m+1}\circ \ga$ 
is simply the restriction of $\rho_{2m+1}$ to $\cl^0_{2m+1}$. But this 
representation is irreducible (see \cite{LM}). The representation 
$\rho_{2m+1}\circ \ga$ is then an irreducible representation of $\cl_{2m}$ of 
dimension $\mathrm{dim}\gS_{2m+1}=2^{[\frac{2m+1}{2}]}=2^m$, as $\rho_{2m}$. 
Now, (\ref{cl2m}) implies that such a representation is unique, up to an 
isomorphism. So $\rho_{2m}\cong \rho_{2m+1}\circ\ga$ and we can conclude that
\begin{equation}\label{dimpair}
S\cong\gS M. 
\end{equation}

Let $\go_{2m}=i^m e_1\cdot\dots\cdot e_{2m}$ be the complex volume form in even 
dimension. An easy calculation shows that $\ga(\go)=\go$. The decomposition of 
$\gS M$ into positive and negative parts is preserved under the isomorphism 
(\ref{dimpair}) and we have $$S=S^+\oplus S^-$$ where $$S^\pm=\set{\psi \in S 
\mid i\nu\cdot\psi=\pm\psi}\cong\gS M^\pm.$$ 

Indeed, because we choose $\rho_{2m+1}$ as the irreducible representation of 
$\cl_{2m+1}$ for which the complex volume form $\go_{2m+1}=i^{m+1} 
e_1\cdot\ldots \cdot e_{2m}\cdot\nu$ acts as the identity on $\gS_{2m+1}$, one 
has, for $\psi \in S$: $$i\nu \cdot\psi=i\nu\cdot\go_{2m+1}\cdot\psi=i^m i^2 
\nu\cdot e_1\cdot \ldots\cdot e_{2m}\cdot\nu\cdot\psi=\go_{2m}\cdot\psi.$$ 

Assume now that $n=2m+1$ is odd. Recall the following isomorphism: 
\begin{equation}\label{climp}
\cl_{2m+1}=\cx(2^m)\oplus\cx(2^m). 
\end{equation}
As mentioned above, $\rho_{2m+1}$ corresponds to the irreducible representation 
of $\cl_{2m+1}$ for which the action of the complex volume form $\go_{2m+1}$ is 
the identity. Because $n+1=2m+2$ is even, $\gS N$ decomposes into positive and 
negative parts, $$\gS N^\pm=\spin N\times_{\rho^\pm_{2m+2}} \gS^\pm_{2m+2}.$$ 

\noindent If $e_k$ is a basis vector tangent to $M$, then 
\begin{eqnarray}
\ga(e_k)\cdot\go_{2m+2}&=&i^{m+1}e_k\cdot\nu\cdot e_1\cdot\ldots\cdot 
e_{2m+1}\cdot \nu\nonumber\\ &=&i^{m+1}(-1)^{2m+2}(-1)^{2m+2}e_1\cdot\ldots\cdot 
e_{2m+1}\cdot\nu\cdot e_k\cdot\nu\nonumber\\&=&\go_{2m+2}\,\ga(e_k).\nonumber 
\end{eqnarray}

So $\rho_{2m+2}\circ \ga$ preserves the decomposition of $\gS N$, and 
$$S=S^+\oplus S^-$$ with $$S^\pm=\spin 
M\times_{\rho^\pm_{2m+2}\circ\ga}\gS_{2m+2}^\pm,$$ and where $\go_{2m+2}$ acts as 
$\pm\mathrm{Id}$ on $S^\pm$. 

Moreover,$$\ga(\go_{2m+1})=i^{m+1} 
(e_1\cdot\nu)\cdot\ldots\cdot(e_{2m+1}\cdot\nu)=i^{m+1}e_1\cdot\ldots\cdot 
e_{2m+1}\cdot\nu=\go_{2m+2},$$ and then $\rho_{2m+1}$ and $\rho^+_{2m+2}\circ 
\ga$ are both irreducible representations of $\cl_{2m+1}$ 
of the same dimension, such that $\rho_{2m+1}(\go_{2m+1})$ and $\rho^+_{2m+2}\circ 
\ga(\go_{2m+1})$ are respectively the identity on $\gS_{2m+1}$ and 
$\gS_{2m+2}^+$. Because such a representation is unique up to an isomorphism, we 
deduce that $\rho_{2m+1}\cong\rho^+_{2m+2}\circ\ga$ and 
\begin{equation}\label{dimimp}
S^+\cong \gS M. 
\end{equation}

Thus we have shown the following
\begin{prop}\label{ident}  
If $n$ is even (resp. odd), there exists an identification of the hypersurface 
spinor bundle $S$ (resp. $S^+$) with the spinor bundle $\gS M$ which sends every 
spinor $\gp \in S$ (resp. $S^+$) to the spinor denoted by $\gp^* \in \gS M$. 
Moreover, with respect to this identification, Clifford multiplication by a 
vector field $X$, tangent to $M$, is given by 
$$X\cdot\gp^*=(X\cdot\nu\cdot\gp)^*.$$ 
\end{prop} 

\subsection{The Spinorial Gauss Formula and the Hypersurface Dirac Operator}

Let $\wit{\gn}$ be the Levi-Civita connection of $(N^{n+1},\wit{g})$, and $\gn$ 
that of $(M^n,g)$. Let $(e_1,\ldots,e_n,e_{n+1}=\nu)$ be a local orthonormal 
basis for $T M$, then the Gauss formula says that for $1\leq i,j \leq n$, 
\begin{equation}\label{gauss}
\wit{\gn}_i e_j=\gn_i e_j + h_{ij}\nu\; ,
\end{equation} 
where $h_{ij}$ are the coefficients of the second fundamental form of the 
hypersurface $M$. We are going to relate the associated connections on the 
corresponding spinor bundles. For this, consider $\phi\in\gG (\gS N)$ and 
$\gp=\phi_{\mid M}\in \gG(S)$ its restriction to $M$. Recall now that locally, 
for $X\in \gG(T M)$, 
\begin{equation}\label{connec1}
\wit{\gn}_X\phi=X(\phi)+\demi\sum_{1\leq i<j \leq n+1} \wit{g}( \wit{\gn}_X 
e_i,e_j)e_i\cdot e_j \cdot \phi 
\end{equation}
and
\begin{eqnarray}\label{connec2}
\gn_X\gp&=&X(\gp)+\demi\sum_{1\leq i<j \leq n} g(\gn_X e_i,e_j)e_i\cdot\nu \cdot 
e_j \cdot\nu \cdot \gp\nonumber\nonumber\\ &=&X(\gp)+\demi\sum_{1\leq i<j \leq n} g(\gn_X 
e_i,e_j)e_i\cdot e_j \cdot \gp.\nonumber
\end{eqnarray}

Therefore, by restricting both sides of equation (\ref{connec1}) to $M$, and 
using the fact that $X(\phi)_{\mid M}=X(\phi_{\mid M})$ for $X$ tangent to $M$, 
the Gauss formula (\ref{gauss}) yields, for $1\leq k \leq n$, 
\begin{eqnarray}
(\wit{\gn}_k\phi)_{\mid M}&=&e_k(\gp)+\demi\sum_{1\leq i<j \leq n} \wit{g}( 
\gn_k e_i + h_{ki}\nu,e_j)e_i\cdot e_j \cdot \gp \nonumber\\&&+ \demi\sum_{1\leq 
i \leq n}\wit{g}( \gn_k e_i + h_{ki}\nu,\nu)e_i\cdot \nu \cdot \gp 
\nonumber\\&=&e_k(\gp)+\demi\sum_{1\leq i<j \leq n} g(\gn_k e_i,e_j)e_i\cdot e_j 
\cdot \gp \nonumber\\&&+ \demi\sum_{1\leq i \leq n}h_{ki}e_i\cdot \nu \cdot \gp 
\nonumber\\&=&\gn_k\gp+\demi\sum_{1\leq i \leq n}h_{ki}e_i\cdot \nu \cdot \gp 
\nonumber \end{eqnarray} 

Once again, from equation (\ref{connec1}), it makes sense to write 
$(\wit{\gn}_X\phi)_{\mid M}=\wit{\gn}_X\gp$ when $X$ is tangent to $M$, and 
hence we proved the spinorial Gauss formula:
\begin{equation}\label{spingauss}
\forall \gp\in\gG(S)\;,\;\forall X\in\gG(T M),\qquad 
\wit{\gn}_X\gp=\gn_X\gp+\demi h(X)\cdot\nu\cdot\gp. 
\end{equation}
(Here $h$ is seen as an endomorphism of the tangent bundle.)

It is known (see \cite{LM}) that there exists a positive definite Hermitian 
metric $<.,.>$ on $\gS N$ such that, if $\tau$ is a $k$-form on $N$, 
\begin{equation}\label{pdscal1}
<\tau\cdot\phi,\psi>=(-1)^{\frac{k(k+1)}{2}}<\phi,\tau\cdot\psi>\quad,\quad 
\forall \phi,\psi\in\gG(\gS N). 
\end{equation}
If we denote $(.,.)$ its real part, we have
\begin{equation}\label{pdscal2}
(\wit{X}\cdot\phi,\wit{Y}\cdot\phi)=\wit{g}(\wit{X},\wit{Y})(\phi,\phi)\quad 
,\quad (\wit{X}\cdot\phi,\phi)=0 \quad ,\quad \forall\wit{X},\wit{Y}\in \gG(T 
N). 
\end{equation}

We simply restrict $(.,.)$ to $M$ to get a globally defined metric on $S$. Now, 
because $\wit{\gn}$ is compatible with $(.,.)$, i.e. 
$$X(\gp,\psi)=(\wit{\gn}_X\gp,\psi)+(\gp,\wit{\gn}_X\psi)\quad ,\quad \forall 
\gp,\psi\in \gG(S),\;\forall X \in \gG(T M)$$ formula (\ref{spingauss}) 
easily implies that $\gn$ is also compatible with the metric. We remark that 
equation (\ref{connec2}) implies that with respect to the identification of 
proposition \ref{ident}, we have 
\begin{equation}\label{iden2}
 (\gn\phi)^*=\gn\phi^*
\end{equation}
This leads to the metric $(.\, ,.)_{\gS M}$ on the intrinsic spinor bundle, with 
the same properties as $(.\, ,.)$, and hence the two bundles are isometric. 

Because Clifford multiplication of a spinor by a vector tangent to $M$ is given 
by (\ref{mulcli}), if $n$ is odd, $S^+$ is stable by $\gn$ and by Clifford 
multiplication. So the classical Dirac operator is simply defined on $S$ for $n$ 
even (resp. $S^+$ for $n$ odd) by $$D=\sum_{i=1}^n e_i\cdot\nu\cdot\gn_i.$$ 

Now we define the hypersurface Dirac operator on $\gG(S)$ as $$D_H=\sum_{i=1}^n 
e_i\cdot\nu\cdot\wit{\gn}_i.$$ 

This definition is motivated by the following fact. Let $$\wit{D}=\sum_{i=1}^n 
e_i \cdot \wit{\gn}_i$$ be the hypersurface Dirac operator defined by Witten 
(see \cite{Wi}, \cite{PT}) to prove the positive energy conjecture in general 
relativity. Then $\wit{D}$ is not formally self-adjoint with respect to the metric 
$(.,.)$. Indeed, it is proved in \cite{HZ1} that $$D_H^2=\wit{D}^{\,*}\wit{D},$$ 
where $\wit{D}^{\,*}$ is the formal adjoint of $\wit{D}$ w.r.t. $(.\, ,.)$. 

From formula (\ref{spingauss}), we see that for $n$ even (resp. odd), we have 
the following relations on $\gG(S)$ (resp. $\gG(S^+)$) : 
\begin{eqnarray}
  D_H&=&\sum_ie_i\cdot\nu\cdot\gn_i+\sum_i e_i\cdot\nu\cdot\frac{h(e_i)}{2}\cdot\nu\cdot\nonumber\\
  &=&D+\sum_{i,j}\frac{h_{ij}}{2}e_i\cdot e_j\cdot\nonumber\\
  &=&D+\sum_{i,j}\frac{h_{ij}}{4}(e_i\cdot e_j+ e_j\cdot e_i)\cdot \nonumber\\
  &=&D-\sum_{i,j}\frac{h_{ij}}{2}\delta_{ij}\nonumber
\end{eqnarray}
and hence, if $H=\sum_i h_{ii}$ is the mean curvature of the hypersurface, we 
have
\begin{equation}\label{DH}
  D_H=D-\frac{H}{2}
\end{equation}
In the following, we will not distinguish the cases where $n$ is even or odd. In 
fact, if $n$ is odd, $D_H$ preserves the decomposition of $S$ into positive and 
negative spinors, as well as Clifford multiplication (recall (\ref{mulcli})), 
$\wit{\gn}$ and $\gn$. Indeed, if $\phi\in \gG(S)$ is an eigenspinor of $D_H$ 
with eigenvalue $\gl$, it is the same for $\phi^+$, its positive part. So we 
only consider positive spinors. The notation becomes easier with this 
convention. 

Now, it is easy to see from equation (\ref{DH}) that $D_H$ is formally self-adjoint with 
respect to the metric $(.,.)$ (see \cite{HZ1}). Finally, recall the well-known 
Schr\"{o}dinger-Lichnerowicz formula on $\gG(\gS M)$ which by the previous 
identification is also true on $\gG(S)$: 
\begin{equation}\label{lichs}
D^2=\gn^*\gn+\frac{R}{4},
\end{equation}
$R$ being the scalar curvature of $M$ and $\gn^*$ the formal adjoint of $\gn$ 
with respect to the metric $(.,.)$.    
 
\section{Proof of the Theorem \ref{mthm}}

Now we give an estimate for the eigenvalues of $D_H$ in terms of the 
energy-momentum tensor (see \cite{Hi3}). For any spinor field $\gp\in \gG(S)$, 
we define the associated energy-momentum 2-tensor $Q^{\,\gp}$ on the complement 
of its zero set, by 
\begin{eqnarray}
Q^{\,\gp}_{ij}=\frac{1}{2} (e_i \cdot \nu \cdot\gn _j \gp + e_j \cdot\nu 
\cdot\gn_i \gp, \gp/|\gp |^2). \label{gq} 
\end{eqnarray}

\begin{rem}
This definition corresponds to the one given in \cite{Hi3} if we note that with 
respect to the identification of $S$ with $\gS M$ of proposition \ref{ident}, 
$$Q^{\,\gp}_{ij}=\frac{1}{2} (e_i\cdot\gn _j \gp^* + e_j\cdot\gn_i \gp^*, 
\gp^*/|\gp^* |^2)_{\gS M}.$$ 
\end{rem}
If $\gp$ is an eigenspinor for $D_H$, $Q^{\,\gp}$ is well defined in the sense 
of distribution. For any real functions $p$ and $q$, consider the modified 
covariant derivative defined on $S$ by 
\begin{eqnarray}
\gn _i ^Q=\gn _i + \Big(p\,\frac{H}{2}+q\,\gl\Big )e_i \cdot \nu \cdot +\,\sum_j 
Q^{\,\gp}_{ij}\;e_j\cdot\nu\cdot \,.\label{gct} 
\end{eqnarray}

\begin{rem}
This connection is well defined on $S^+$ when $n$ is odd. 
\end{rem}
Using (\ref{pdscal2}), we have
\begin{eqnarray*}
|\gn^Q \gp |^2 &=& |\gn \gp |^2+n\;\Big(p\;\frac{H}{2}+q\; \gl\Big )^2\;|\gp 
|^2\\ &&+\sum_{i,j,k} Q^{\,\gp}_{ij}\, Q^{\,\gp}_{ik} (e_j \cdot\nu\cdot 
\gp,\;e_k\cdot \nu \cdot \gp)\\ && +2\;\Big(p\;\frac{H}{2}+q\; \gl\Big )\;\sum 
_i (\gn _i \gp,\; e_i \cdot \nu \cdot\gp)\\ &&+2\sum_{i,j} Q^{\,\gp}_{ij} \;( 
\gn _i \gp,\; e_j\cdot\nu \cdot \gp)\\ && +2\Big(p\;\frac{H}{2}+q\;\gl\Big 
)\;\sum _{i,j} Q^{\,\gp}_{ij}\; (e_i \cdot \nu \cdot \gp , 
e_j\cdot\nu\cdot\gp)\\ &=& |\gn \gp|^2+n\;\Big(p\;\frac{H}{2}+q\; \gl\Big 
)^2\;|\gp |^2 +|Q^{\,\gp}| ^2 \;|\gp |^2\\ && -2 
\;\Big(p\;\frac{H}{2}+q\;\gl\Big )\;  (D \gp,\; \gp)-2|Q ^\gp| ^2 \;|\gp |^2\\ 
&& +2\Big(p\;\frac{H}{2}+q\; \gl\Big )\; Tr(Q^{\,\gp})\;|\gp |^2 
\end{eqnarray*}
but $$Tr(Q^{\,\gp})\;|\gp |^2=(D \gp,\gp)$$ hence 
\begin{equation}\label{Qformula}
|\gn^Q \gp |^2 = |\gn \gp|^2+n\;\Big(p\;\frac{H}{2}+q\; \lambda\Big )^2\;|\gp 
|^2 -|Q^{\,\gp}| ^2 \;|\gp |^2. 
\end{equation}

Now, since $D_H=D-\frac{H}{2}$, the Schr\"odinger-Lichnerowicz formula 
(\ref{lichs}) on $\gG(S)$ gives 
\begin{eqnarray}
\int_M |\gn\gp|^2 v_g 
&=&\int_M\Big(|D\gp|^2-\frac{R}{4}|\gp|^2\Big)v_g\nonumber\\ 
&=&\int_M\Big((\gl+\frac{H}{2})^2-\frac{R}{4}\Big)|\gp|^2v_g\label{lich} 
\end{eqnarray}

Therefore (\ref{Qformula}) and (\ref{lich}) imply 
\begin{eqnarray}\label{Qint}
\int_M |\gn^Q \gp |^2v_g&=& \int_M \Big( (1+nq^2)\gl^2-\frac{R}{4}-|Q^{\,\gp}| ^2 \Big)|\gp |^2 v_g\nonumber\\
&&+\int_M\Big((1+np^2)\frac{H^2}{4}+(1+npq)H\gl\Big)|\gp |^2 v_g.  
\end{eqnarray}

Now, assume that $q$ has no zeros so that we can choose $p=-\frac{1}{nq}$. Then 
(\ref{Qint}) becomes 
\begin{equation}
  \label{Qpres}
 \int_M |\gn^Q \gp |^2v_g= \int_M
(1+nq^2)\Big[\gl^2-\frac{1}{4}\Big(\frac{R+4|Q^{\,\gp}| 
^2}{(1+nq^2)}-\frac{H^2}{nq^2}\Big)\Big] |\gp |^2 v_g 
\end{equation}

If $R+4|Q^{\,\gp}| ^2>H^2>0$, we can take 
\begin{equation}
  \label{nq2}
  nq^2=\frac{|H|}{\sqrt{R+4|Q^{\,\gp}| ^2}-|H|}
\end{equation}

Then equation (\ref{Qpres}) becomes 
\begin{equation}
  \label{Qoui}
\int_M |\gn^Q \gp |^2v_g= \int_M (1+nq^2)\Big[\gl^2-\frac{1}{4}\left( \sqrt {R+4 
|Q^{\,\gp}| ^2} -|H| \right) ^2 \Big] |\gp |^2 v_g 
\end{equation}

Because the left hand side of this equation is positive and $\gl$ is a constant, 
we get 
\begin{equation}\label{theo}
\gl^2 \geq \frac{1}{4}\; \inf _M \left( \sqrt {R+4 |Q^{\,\gp}| ^2} -|H| \right) ^2. 
\end{equation}

\begin{rem}
If $M$ is minimal, i.e. $H=0$, we can choose $q\equiv 0$ in (\ref{gct}) so that 
(\ref{theo}) specializes to the inequality of Theorem A in \cite{Hi3}. 
\end{rem} 

\begin{rem}
Note that our definition of the energy-momentum tensor $Q^{\,\gp}$ coincides 
with that in \cite{Hi3}. The definition used in \cite{Z1} and \cite{HZ1} gives a 
factor $\frac{n}{n-1}$ in front of $R+4 |Q^{\,\gp}| ^2$ in inequality 
(\ref{theo}) but in this case, $Q^{\,\gp}$ has no canonical intrinsic meaning. 
\end{rem}

\section{Limiting cases}

First recall the inequality proved by X. Zhang :

\begin{thm}[\cite{Z1}]\label{zhang}
Let $M^n \subset N^{n+1}$ be a compact hypersurface of a riemannian spin 
manifold $(N,\wit{g})$. Assume that $n\geq 2$ and $n R>(n-1) H^2 > 0$. Then if 
$\gl$ is any eigenvalue of the hypersurface Dirac operator $D _H$, one has 
\begin{equation}\label{zbin}
\gl^2 \geq \frac{1}{4}\; \inf _M \left( \sqrt {\frac{n}{n-1}R} -|H| \right) ^2. 
\end{equation} 
\end{thm}

As in the proof of Theorem \ref{mthm}, the proof of Theorem \ref{zhang} is based 
on the use of the modified connection 
\begin{equation}\label{zgn}
\gn _i ^\gl=\gn _i + \Big(p\,\frac{H}{2}+q\,\gl\Big )e_i \cdot \nu \cdot.
\end{equation} 
Here, $p$ and $q$ are related by 
\begin{equation}\label{defp}
  p=\frac{1-q}{1-n q}
\end{equation}
and
\begin{equation}\label{defqq}
  q=\frac{1}{n}\Big(1-\sqrt{\frac{(n-1)|H|}{\sqrt{\frac{n}{n-1}R} -|H|}}\Big)
\end{equation}
or, in other terms,
\begin{equation}\label{defq}
  (1-n q)^2= \frac{(n-1)|H|}{\sqrt{\frac{n}{n-1}R} -|H|}.
\end{equation}
Equality holds in (\ref{zbin}) for an eigenspinor $\gp$ of $D_H$ with eigenvalue 
$\gl$ if and only if $\sqrt{\frac{n}{n-1}R} -|H|$ is constant and 
$\gn^\gl\gp\equiv 0$.  But, with respect to the identification of proposition 
\ref{ident}, and by (\ref{iden2}), $\gn^\gl\gp\equiv 0$ is equivalent to 
\begin{equation}\label{killing}
\forall i=1,\ldots, n,\qquad\gn_i\gp^*=-\Big(p\,\frac{H}{2}+q\,\gl\Big )e_i 
\cdot\gp^* 
\end{equation}
It is known (see \cite{Hi1}) that if such a section exists on $\gS M$, then 
$p\,\frac{H}{2}+q\,\gl$ has to be constant (say $\frac{\gl_1}{n}$ for instance) 
and that in this case $M$ is Einstein and $R=4\frac{n-1}{n}\gl_1^2$. So $\gp$ is 
a Killing spinor and we are in the limiting case of Friedrich's inequality 
\cite{Fr}. Moreover, since $\sqrt{\frac{n}{n-1}R} -|H|$ is constant, $H$ has to 
be constant. 

Therefore, since $D\gp=\gl_1\gp$ and $\gl_1=\frac{\mathrm{sign}(\gl_1)}{2}\sqrt 
{\frac{n}{n-1}R}$, the following equation must be satisfied (recall that 
$D_H=D-\frac{H}{2}$) 
\begin{equation}\label{temp}
\gl=\frac{\mathrm{sign}(\gl_1)}{2}\sqrt 
{\frac{n}{n-1}R}-\frac{H}{2}=\frac{\mathrm{sign}(\gl_1)}{2}\sqrt 
{\frac{n}{n-1}R}-{\mathrm{sign}(H)}\frac{|H|}{2} 
\end{equation} 
But equality case gives 
\begin{equation}\label{temp1}
\gl=\frac{\mathrm{sign}(\gl)}{2}\Big(\sqrt {\frac{n}{n-1}R} -|H|\Big),   
\end{equation}
So (\ref{temp}) and (\ref{temp1}) imply that
\begin{equation}\label{sign}
  \mathrm{sign}(\gl)=\mathrm{sign}(\gl_1)=\mathrm{sign}(H)
\end{equation}
On the other hand, an easy calculation leads to
\begin{eqnarray}\label{sign1}
p\,\frac{H}{2}+q\,\gl&=&\frac{\mathrm{sign}(\gl)}{2n} 
\sqrt{\frac{n}{n-1}R}\nonumber\\ &&+\frac{(\mathrm{sign}(H)- 
\mathrm{sign}(\gl))}{2n}\Big(1+\sqrt{(n-1)(\sqrt{\frac{n}{n-1}R} -|H|)}\Big) 
\nonumber\\&=&\frac{\mathrm{sign}(\gl_1)}{2n} \sqrt{\frac{n}{n-1}R}\nonumber 
\end{eqnarray}
and we recover the already known fact that 
$p\,\frac{H}{2}+q\,\gl=\frac{\gl_1}{n}$. 

Indeed, (\ref{sign}) can be trivially observed because in the equality case, 
both $R$ and $H$ are constant, so we can think of the spectrum of $D_H$ as the 
shifting of the spectrum of $D$ by $-\frac{H}{2}$. Then the condition $n R>(n-1) 
H^2 > 0$ in Theorem \ref{zhang} simply implies that the lowest eigenvalue of 
$D_H$ (in the sense of its absolute value) must have the sign of $H$. In 
particular, when $n$ is even, it shows how we lose the symmetry of the spectrum when passing from 
$D$ to $D_H$ (compare with the case where $H=0$). 
 
Now we discuss the case of Theorem \ref{mthm}. The limiting case of inequality 
(\ref{bin}) holds for an eigenspinor $\gp$ of $D_H$ with eigenvalue $\gl$ if and 
only if $\gn^Q\gp \equiv 0$. First note that this implies that $|\gp|$ is 
constant. Then, with respect to the identification of proposition \ref{ident}, 
and by (\ref{iden2}), $\gn^Q\gp\equiv 0$ is equivalent to 
\begin{equation}\label{tk}
\gn_i\gp^*=-\Big(p\,\frac{H}{2}+q\,\gl\Big )e_i \cdot\gp^*-\sum_j 
Q^{\,\gp}_{ij}\;e_j\cdot\gp^* 
\end{equation}
Let $f=p\,\frac{H}{2}+q\,\gl$, then equation (\ref{tk}) can be written as 
\begin{equation}\label{tk2}
\gn_i\gp^*=-\sum_j (Q^{\,\gp}_{ij}+f\delta_{ij})\;e_j\cdot\gp^*  
\end{equation}
  
Now let $T_{ij}=Q^{\,\gp}_{ij}+f\delta_{ij}$, taking Clifford multiplication by 
$e_k$ on both sides of equation ($\ref{tk2}$), yields 
$$e_k\cdot\gn_i\gp^*=-\sum_j T_{ij}\, e_k\cdot e_j\cdot \gp^* ,$$ which gives 
$$(e_k\cdot\gn_i\gp^*,\gp^*)_{\gS M}=-\sum_j T_{ij}(e_k\cdot e_j\cdot \gp^*, 
\gp^*)_{\gS M},$$ and, because $(e_k\cdot e_j\cdot \gp^*, \gp^*)_{\gS M}=0$ 
unless $j=k$ and $T_{ij}$ is symmetric, we proved 
$$\demi(e_i\cdot\gn_k\gp^*+e_k\cdot\gn_i\gp^*,\gp^*/|\gp^*|^2)_{\gS M}=T_{ik}.$$ 
Hence $$T_{ik}=Q^{\,\gp}_{ik}$$ and we can conclude that $f=0$. Equation 
(\ref{tk}) reduces to 
\begin{equation}\label{tk3}
\gn_i\gp^*=-\sum_j Q^{\,\gp}_{ij}\;e_j\cdot\gp^*. 
\end{equation}
Such field equations have been studied, as well as their integrability 
conditions, by T. Friedrich and E. C. Kim in \cite{FK1}. Note that they allow a 
nice formulation of the theory of immersed surfaces in the euclidean 3-space 
(see \cite{Fr}). We will call an EM-spinor (for Energy-Momentum spinor) a non 
trivial spinor field satisfying (\ref{tk3}). If it is an eigenspinor for the 
Dirac operator, which is equivalent to the fact that $\mathrm{tr}\,Q^{\,\gp}$ is 
constant, it is called T-Killing spinor (see \cite{FK}). In fact, a T-Killing 
spinor is exactly a spinor field satisfying the limiting case in Hijazi's 
inequality \cite{Hi3}. 

Now we have (see \cite{Hi3} or Lemma 4.1(iii) of \cite{FK1})
\begin{equation}\label{Qtr}
  (\mathrm{tr}\,Q^{\,\gp})^2=\frac{R}{4}+ |Q^{\,\gp}| ^2
\end{equation}
So (\ref{tk3}) implies that $$D\gp=F\gp$$ where $F^2=\frac{R}{4}+ |Q^{\,\gp}| 
^2$. Whereas equality case in (\ref{bin}) gives $\sqrt {R+4 |Q^{\,\gp}| ^2} 
-|H|$ is constant, we can't conclude here that $\frac{R}{4}+ |Q^{\,\gp}| ^2$ and 
$H$ are constant as in the case of Zhang's inequality. Nevertheless, we have the 
following 
\begin{cor}
If $H$ is constant, then equality case in (\ref{bin}) holds if and only $\gp$ is 
a T-Killing spinor. 
\end{cor}   
By hypothesis $H$ has constant sign and we can conclude that $\gl$ has the same 
sign. Recall that $p$ and $q$ are related by $$p=-\frac{1}{nq}$$ and 
$$nq^2=\frac{|H|}{\sqrt{R+4|Q^{\,\gp}| ^2}-|H|}.$$ Indeed, an easy calculation 
gives 
\begin{equation}\label{sign3}
0=f=\Big(p\,\frac{H}{2}+q\,\gl\Big )=\frac{(\mathrm{sign}(\gl)- 
\mathrm{sign}(H))}{2\sqrt{n}}\sqrt{|H|(\sqrt{R+4|Q^{\,\gp}| ^2} -|H|)}. 
\end{equation}
Hence $$\mathrm{sign}(\gl)= \mathrm{sign}(H).$$ 

\begin{rem}\label{srem}
Equality case of (\ref{zbin}) is included in that of (\ref{bin}): if we assume 
that $\gp$ is a Killing spinor, then necessarily 
$Q^{\,\gp}_{ij}=\frac{\gl_1}{n}\delta_{ij}$ and so 
$(\mathrm{tr}\,Q^{\,\gp})^2=\gl_1^2=\quart\frac{n}{n-1}R$. Therefore equation 
(\ref{Qtr}) implies $$ 4|Q^{\,\gp}| ^2=\frac{n}{n-1}R-R $$ and we have $$\gl^2= 
\left( \sqrt {\frac{n}{n-1}R} -|H| \right) ^2.$$ 
\end{rem}

\begin{rem}
The previous remark shows that Theorem \ref{mthm} improves Theorem \ref{zhang}. 
In particular, it does not require $R$ to be positive, and the limiting case 
does not imply that $H$ has to be constant. \end{rem}

\section{Proof of the Theorem \ref{sthm}}

Consider a conformal change of the metric $\ovl{g}=e^{2u}\wit{g}$ for any real 
function $u$ on $N$. For simplicity, let $\ovl{N}=(N,\ovl{g})$. The natural 
isometry between $\so N$ and $\so \ovl{N}$ induced by this conformal change of 
the metric lifts to an isometry between the $\spin (n+1)$-principal bundles 
$\spin N$ and $\spin \ovl{N}$, and hence between the two corresponding 
hypersurface spinor bundles $S$ and $\ovl{S}$. If $\gp\in\gG (S)$, denote by 
$\ovl{\gp}\in\gG(\ovl{S})$ its image by this isometry. Let $(.,.)_{\ovl{g}}$ be 
the metric on $\ovl{S}$ naturally defined as described in section 2. Then for 
$\gp$, $\psi$ two sections of $S$, we have 
$$(\gp,\psi)=(\ovl{\gp},\ovl{\psi})_{\ovl{g}}\quad \mathrm{and}\quad 
\ovl{X}\,\bar{\cdot}\,\ovl{\psi}=\ovl{X\cdot\psi} $$  

We will also denote by $\ovl{g}={e^{2u}}_{\mid M}g$ the restriction of $\ovl{g}$ 
to $M$. By conformal covariance of the Dirac operator, we have, for $\gp \in 
 \gG(S)$, (see \cite{HZ1})
\begin{equation}\label{covconf}
\overline D \Big(\; e ^{-\frac{(n-1)}{2} u} \;\overline{\gp}\; \Big)=e 
^{-\frac{(n+1)}{2} u} \;\overline {D \gp},  
\end{equation}
where $\overline D$ stands for the Dirac operator w.r.t. to $\ovl{g}$. On the 
other hand 
\begin{equation}\label{Hbar}
\overline  H = e ^{-u} \,\Big( H+ n\; du (\nu) \Big). 
\end{equation}
Therefore, if $\overline{D}_{\ovl{H}}$ stands for the hypersurface Dirac 
operator w.r.t. to $\ovl{g}$, equations (\ref{covconf}) and (\ref{Hbar}) imply 
that, 
\begin{equation}\label{DHbar}
\ovl{D}_{\ovl{H}}\; \Big(\; e ^{-\frac{(n-1)}{2} u}\; \ovl{\gp}\; \Big) = e 
^{-\frac{(n+1)}{2} u }\; \Big(\; {\overline {D _H \gp}}-\frac{n}{2}\;du(\nu)\; 
\ovl{\gp} \;\Big). 
\end{equation}
\begin{rem}\label{reg}
We see that if $du (\nu)_{\mid M} =0$, $D_H$ is a conformal invariant operator. 
In this case, technics used in (\cite{Hi1}) can be applied for the eigenvalues 
of $D_H$. Indeed, such a conformal change of metric can be viewed as a intrinsic 
conformal change of the metric on $M$, when we omit the ambient space $N$ (See 
section 7). 
\end{rem}
From now on, we will only consider conformal changes of the metric 
$\ovl{g}=e^{2u}\wit{g}$ with $du (\nu)=0$ on $M$. They will be called regular 
conformal changes of metric as in \cite{HZ1}. 

For $\gp \in \gG(S)$ an eigenspinor of $D_H$ with eigenvalue $\gl$, let 
$\ovl{\psi}:= e ^{-\frac{n-1}{2} u}\, \overline \gp$. Then (\ref{DHbar}) gives
\begin{equation}\label{confeig}
\overline{D}_{\ovl{H}} \;\overline \psi = \lambda\; e ^{-u}\; \overline  
\psi\, . 
\end{equation}

Recall that
\begin{eqnarray}\label{gnconf}
\overline \gn _i \overline \gp = \overline {\gn _i \gp} -\frac{1}{2} \;\overline 
{e_i \cdot \mathrm{d} u \cdot \gp} -\frac{1}{2} e_i(u) \; \overline \gp, 
\end{eqnarray}
and $\ovl{e_i}=e^{-u}e_i.$ Now, as in \cite{Hi1}, it is straightforward to get  
\begin{eqnarray}\label{QQQ}
\overline Q^{\,\overline \psi}_{\bar i\, \bar j} &=&\frac{1}{2}\;  \Big 
(\overline {e_i}\; \bar \cdot\; \overline {\nu}\; \bar \cdot \;\overline \gn 
_{\overline e _j} \;\overline \psi \;+\; \overline {e _j}\; \bar \cdot\;  
\overline {\nu}\; \bar \cdot\; \overline \gn _{\overline e _i} \;\overline 
\psi\;, \;\overline \psi/|\overline \psi |^2 _{\overline g} \Big )_{\ovl{g}} 
\nonumber\\ &=& \frac{1}{2}\; e ^{-u} \; \Big (\overline {e _i} \;\bar \cdot 
\;\overline {\nu} \;\bar \cdot \;\overline \gn _{e _j} \;\overline \gp\; + 
\;\overline {e _j} \;\bar \cdot \;\overline {\nu} \;\bar \cdot \;\overline \gn 
_{e _i} \;\overline \gp\;, \;\overline \gp/|\overline \gp |^2 _{\overline g} 
\Big )_{\ovl{g}}\nonumber\\ &=& \frac{1}{2}\; e ^{-u} \; \Big (e _i\cdot\nu 
\cdot\gn _{e _j}  \gp\; + e_j\cdot\nu\cdot \gn _{e _i} \gp\;, \;\gp/|\gp |^2 
\Big )\nonumber\\ 
 &=& e ^{-u}Q ^{\,\gp}_{ ij}. 
\end{eqnarray}
Hence,
\begin{equation}\label{QQ}
|\overline Q^{\,\overline \psi}|^2=e^{-2u}|Q ^{\,\gp}|^2
\end{equation}
Equation (\ref{Qpres}), which is also true on $N$, applied to $\ovl{\psi}$ 
yields 
\begin{equation}
 \int_M |{\ovl{\gn}}^Q\, \ovl{\psi}|^2 v_{\ovl{g}}= \int_M
(1+nq^2)\Big[\gl^2e^{-2u}-\frac{1}{4}\Big(\frac{\ovl{R}+4|\ovl{Q}^{\,\ovl{\psi}}| 
^2}{(1+nq^2)}-\frac{\ovl{H}^2}{nq^2}\Big)\Big] |\ovl{\psi} |^2 v_{\ovl{g}} 
\end{equation}
which, because of (\ref{Hbar}) and (\ref{QQ}) gives
\begin{equation}\label{Qbarpres}
 \int_M |{\ovl{\gn}}^Q\, \ovl{\psi}|^2 v_{\ovl{g}}= \int_M
(1+nq^2)\Big[\gl^2-\frac{1}{4}\Big(\frac{\ovl{R}e^{2u}+4|Q^{\,\gp}| 
^2}{(1+nq^2)}-\frac{H^2}{nq^2}\Big)\Big]e^{-2u} |\ovl{\psi} |^2 v_{\ovl{g}}. 
\end{equation}

Taking $$nq^2=\frac{|H|}{\sqrt{\ovl{R}e^{2u}+4|Q^{\,\gp}| ^2}-|H|}$$  completes 
the proof of Theorem \ref{sthm}. 

\section{General limiting cases}

We now discuss the limiting case in inequality (\ref{ssbin}). Equality holds if 
and only if ${\ovl{\gn}}^Q_{\bar{i}}\, \ovl{\psi}=0$ for $1\leq i\leq n$, which 
can be written as $$0={\ovl{\gn}}_{\bar{i}} \ovl{\psi}+ 
\Big(p\,\frac{\ovl{H}}{2}+qe^{-u}\gl\Big )\ovl{e}_i\, \ovl{\cdot}\, \ovl{\nu}\, 
\ovl{\cdot}\, \ovl{\psi}+\, \sum_j {\ovl{Q}}^{\,\ovl{\psi}}_{ij}\; 
\ovl{e}_j\,\ovl{\cdot}\,\ovl{\nu}\,\ovl{\cdot}\,\ovl{\psi}. $$ Since $\ovl{\psi}:= e 
^{-\frac{n-1}{2} u}\, \overline \gp$, (\ref{gnconf}) and (\ref{QQ}) yield 
\begin{eqnarray} 0 &=& e ^{-\frac{n-1}{2} u}e^{-u}\Big [ \overline {\gn _i \gp} 
-\frac{1}{2} \;\overline {e_i \cdot \mathrm{d} u \cdot \gp} -\frac{n}{2} e_i(u) 
\; \overline \gp \nonumber\\&&+\; \Big(p\,\frac{H}{2}+q\,\gl\Big )\ovl{e_i\, 
\cdot\nu\cdot\gp}+\,\sum_j  {Q}^{\,\gp}_{ij}\; \ovl{e_j\cdot\nu\cdot\gp}\Big ]. 
\end{eqnarray}
With respect to the identification of proposition \ref{ident}, and by 
(\ref{iden2}), this last statement is equivalent to  
\begin{eqnarray}\label{presq} \gn _i \gp^* =\frac{1}{2} \;
e_i \cdot \mathrm{d} u \cdot \gp^* + \frac{n}{2} \mathrm{d} u(e_i)\gp^* - f \; 
e_i\cdot\gp^*-\, \sum_j {Q}^{\,\gp}_{ij}\;e_j\cdot\gp^*. 
\end{eqnarray}
where $f:=p\,\frac{H}{2}+q\,\gl$. As in Section 4, let 
$T_{ij}=Q^{\,\gp}_{ij}+f\delta_{ij}$. It is then straightforward to prove that 
$T_{ij}=Q^{\,\gp}_{ij}$ and so $f=0$. 

Taking the scalar product of (\ref{presq}) with $\gp^*$, it follows 
\begin{eqnarray}\frac{1}{2} e_i(|\gp|^2)&=&(\gn_i\gp^*,\gp^*)_{\gS M}\nonumber\\
&=&\frac{1}{2}( \; e_i \cdot \mathrm{d} u \cdot \gp^*, \gp^*)_{\gS M}+\frac{n}{2} 
\mathrm{d} u(e_i)|\gp|^2\nonumber\\ &=& \frac{(n-1)}{2} \mathrm{d} u(e_i)|\gp|^2 
\nonumber\end{eqnarray}

Therefore,
\begin{equation}\label{du}
 \mathrm{d} u=\frac{\mathrm{d}|\gp|^2}{(n-1)|\gp|^2}. 
\end{equation}

So we proved that equality holds in (\ref{ssbin}) if and only if the eigenspinor 
$\gp$ satisfies $$ \gn _i \gp^* =\frac{1}{2} \; e_i \cdot \mathrm{d} u \cdot 
\gp^* + \frac{n}{2}\, \mathrm{d} u(e_i)\,\gp^* -\sum_j 
{Q}^{\,\gp}_{ij}\;e_j\cdot\gp^*$$ with $u$ satisfying (\ref{du}). Such field 
equations have already been studied, as well as their integrability conditions, 
by T. Friedrich and E.C. Kim \cite{FK1}. We will call them WEM-spinors (for Weak 
Energy-Momentum spinors). If, they satisfy the Einstein-Dirac equation, they are 
called WK-spinors (for Weak Killing spinors). These are exactly the limiting 
case of Hijazi's equality involving conformal change of the metric and the 
energy-momentum tensor \cite{Hi3}, in which case, they are also eigenspinors for 
the classical Dirac operator. 

In our situation, there are not eigenspinors for $D$. As a consequence, even if 
in the limiting case $\sqrt{\ovl{R}e^{2u}+4|Q^{\,\gp}| ^2}-|H|$ has to be 
constant, we can't conclude that both $\sqrt{\ovl{R}e^{2u}+4|Q^{\,\gp}| ^2}$ and 
$H$ are constant. 

Nevertheless, as in the previous section, a simple calculation leads to $$ 
0=f=e^u\frac{(\mathrm{sign}(\gl)- 
\mathrm{sign}(H))}{2\sqrt{n}}\sqrt{|H|(\sqrt{\ovl{R}e^{2u}+4|Q^{\,\gp}| ^2} 
-|H|)}. $$ Hence $$\mathrm{sign}(\gl)= \mathrm{sign}(H).$$    

Now recall the inequality proved by O. Hijazi and X. Zhang : 

\begin{thm}[\cite{HZ1}]\label{hzthm}
Let $M^n \subset N^{n+1}$ be a compact hypersurface of a Riemannian spin 
manifold $(N,\wit{g})$. Assume that $n\geq 2$ and $n \,\ovl{R}\,e^{2u}>(n-1) H^2 
> 0$ for some regular conformal change of the metric $\ovl{g}=e^{2u}\wit{g}$. Then 
if $\gl$ is any eigenvalue of the hypersurface Dirac operator $D _H$, one has 
\begin{equation}\label{hzbin}
\gl^2 \geq \frac{1}{4}\; \inf _M \left( \sqrt {\frac{n}{n-1}\,\ovl{R}\,e^{2u}} 
-|H| \right) ^2. 
\end{equation} 
\end{thm}

As in the proof of Theorem \ref{sthm}, Theorem \ref{hzthm} is obtained by using 
the modified connection defined by (\ref{zgn}), on the manifold 
$(N,\ovl{g}=e^{2u}\wit{g})$. 

As in the beginning of this section, it is then easy to see that equality holds 
in (\ref{hzbin}) if and only if  
\begin{eqnarray} 0 = e ^{-\frac{n-1}{2} u}e^{-u}\Big [ \overline {\gn _i \gp} 
-\frac{1}{2} \;\overline {e_i \cdot \mathrm{d} u \cdot \gp} -\frac{n}{2} e_i(u) 
\; \overline \gp +\; \Big(p\,\frac{H}{2}+q\,\gl\Big )\ovl{e_i\, 
\cdot\nu\cdot\gp}\Big ]. 
\end{eqnarray}
With respect to the identification of proposition \ref{ident}, and by 
(\ref{iden2}), this last statement is equivalent to  
\begin{eqnarray}\label{hhzz}\gn _i \gp^* =\frac{1}{2} \;
e_i \cdot \mathrm{d} u \cdot \gp^* + \frac{n}{2} \mathrm{d} u(e_i)\gp^* - f \; 
e_i\cdot\gp^*. 
\end{eqnarray}
where $f:=p\,\frac{H}{2}+q\,\gl$. As in section 4, let $T_{ij}=f \delta_{ij}$. 
Then it is straightforward to prove that $T_{ij}=Q^{\;\gp}_{ij}$ and that 
spinors fields satisfying the equality case in Theorem \ref{hzthm} are 
particular WEM-spinors. Now, by (\ref{Qtr}) and by (\ref{QQQ}), we see that 
necessarily 
\begin{equation}\label{ff}
  f=\pm\frac{1}{n}\sqrt {\frac{n}{n-1}\,\ovl{R}\,e^{2u}}.
\end{equation}
Hence, solutions of (\ref{hhzz}) correspond exactly to sections verifying the 
limiting case of inequality (5.1) in \cite{Hi1}.
 
 Recall that here, 
$p$ and $q$ are given by 
\begin{equation}
  p=\frac{1-q}{1-n q}
\end{equation}
and 
\begin{equation}
  (1-n q)^2= \frac{(n-1)|H|}{\sqrt{\frac{n}{n-1}\,\ovl{R}\,e^{2u}} -|H|}.
\end{equation}
Therefore we can recover that $\mathrm{sign}(\gl)=\mathrm{sign}(H)$ as made 
previously by computing explicitly $p\frac{H}{2}+q\gl$. In fact, as in Remark 
\ref{srem}, the consequence of equation (\ref{Qtr}) is that Theorem \ref{sthm} 
improves Theorem \ref{hzthm}. 
   
\section{Concluding remark} 

We conclude this paper by observing that all computations previously made could 
be done in an intrinsic way, considering a modified Dirac operator 
$D_f=D-\frac{f}{2}$, and connections on $\gS M$: $$ \gn _i ^\gl=\gn _i + 
\Big(p\,\frac{f}{2}+q\,\gl\Big )e_i \cdot $$ and $$ \gn _i ^Q=\gn _i + 
\Big(p\,\frac{f}{2}+q\,\gl\Big )e_i \cdot +\, Q^{\,\gp}_{ij}\;e_j\cdot\, , $$ 
with an appropriate choice of $p$ and $q$ (simply replace $H$ by $f$), in 
(\ref{defp}) and (\ref{defq}) . 
 
The identification of the spinor bundles of section 2 allows to assert that 
computations will lead to the same results, but in a more general way, and so we 
proved 

\begin{prop}
Let $(M^n,g)$ be a compact Riemannian spin manifold. Assume that $n\geq 2$ and 
$n R>(n-1) f^2 > 0$, with $f:M\rightarrow \rl$ a smooth function. Then for any 
eigenvalue $\gl$ of the Dirac-Schr\"{o}dinger operator $D_f=D-\frac{f}{2}$, one has 
$$ \gl^2 \geq \frac{1}{4}\; \inf _M \left( \sqrt {\frac{n}{n-1}R} -|f| \right) 
^2. $$ Equality holds if and only if $M$ admits a Killing spinor and in this 
case $(M^n,g)$ is Einstein, $f$ constant, and 
$$\mathrm{sign}(\gl)=\mathrm{sign}(f).$$ 
\end{prop}
Similarly, one has
\begin{prop}
Let $(M^n,g)$ be a compact Riemannian spin manifold. Let $\gl$ be any eigenvalue 
of the Dirac-Schr\"{o}dinger operator $D_f=D-\frac{f}{2}$, associated with the 
eigenspinor $\gp$. Assume that $R+4 |Q^{\,\gp}| ^2 > f^2 > 0$, then $$ \gl^2 
\geq \frac{1}{4}\; \inf _M \left( \sqrt {R+4 |Q^{\,\gp}| ^2} -|f| \right) ^2. $$ 
where $Q^{\,\gp}$ the energy-momentum tensor associated with $\gp$. 

If equality holds, $M$ admits an EM-spinor, and in this case, 
$$\mathrm{sign}(\gl)=\mathrm{sign}(f).$$ 
\end{prop}

Now using a conformal change of the metric $g$, (see remark \ref{reg}), we prove 
in the same way 
 
\begin{prop}\label{ssthm}
Let $(M^n,g)$ be a compact Riemannian spin manifold. Let $\gl$ be any eigenvalue 
of the Dirac-Schr\"{o}dinger operator $D_f=D-\frac{f}{2}$, associated with the 
eigenspinor $\gp$.\\ If $ \overline R e ^{2u} +4|Q^{\,\gp}| ^2>f ^2>0$, where 
$\overline R$ is the scalar curvature of $M$ for a conformal metric $\overline g 
= e ^{2u} g$, then 
\begin{eqnarray*}
\gl^2 \geq \frac{1}{4}\; \inf _M \left( \sqrt {\overline R\; e ^{2u}+4 
|Q^{\,\gp}| ^2} -|f| \right) ^2. 
\end{eqnarray*}

If equality holds, $M$ admits an WEM-spinor, and in this case, the 
function $u$ is uniquely defined up to a constant by  
$$u=\frac{\ln(|\gp|^2)}{(n-1)}.$$ Moreover,  
$$\mathrm{sign}(\gl)=\mathrm{sign}(f).$$
\end{prop} 

\providecommand{\bysame}{\leavevmode\hbox to3em{\hrulefill}\thinspace}

\end{document}